\documentclass[11pt,a4paper]{article}
\usepackage{amstex}
\usepackage{delarray} 
\usepackage[english,french]{babel}
\newtheorem{defi}{Definition }[section]

\newtheorem{theo}{Theorem }[section]

\newtheorem{cor}{Corollary }[section]
\newtheorem{pro}{Proposition }[section]

\newtheorem{rk}{Remark}[section]
\newtheorem{stp}{Step}[theo]
\newcommand{\C}{{\Bbb C}}
\newcommand{\Z}{{\Bbb Z}}

\newcommand{\N}{{\Bbb N}}
\newcommand{\uq}{\mathrm{u}_q^+}
\newcommand{\sld}{\frak{sl}_2 \left(\Bbb C\right)}

\title{Monoidal structure of the category of $\uq$-modules}
\author{El\'\i sabet Gunnlaugsd\'ottir}
\date{}
\begin{document}
\selectlanguage{english}
\maketitle \footnote{2000 Mathematical subject classification,
primary 20G42 ; secondary 18D10.}

\section{Introduction.}
We consider the half-quantum group $\uq(\sld)$ at a root of unity which order is not $4$.
This non quasi-cocommutative 
Hopf algebra is the upper triangular sub-Hopf algebra of $
\mathrm{u}_q(\sld)$, quotient of the quantized enveloping algebra at a
root of unity $q$ (see \cite{5}). 
Half quantum groups provide
universal $R$-matrices through the Drinfeld double and hence solutions
to the Yang-Baxter equation. Furthermore they appear of interest in
knot theory and 3-manifold invariants. For a simple Lie algebra
$\frak G$, a presentation of $\uq(\frak G)$ by quiver and relations has
been  established by Cibils in \cite{3}, showing that only $\uq$ is of
finite representation type, the others being of tame or wild representation
type. 

In order to study more deeply the representation theory of $\uq$, we
consider the particular family of indecomposable modules on $\uq$
which are $\mathrm{u}_q$-modules as well. We  call them ``extendable
modules'' .  They form a subring of the Grothendieck ring of $\uq$, and
their study leads to a
Clebsch-Gordan-like formula for the decomposition of the tensor
product, taken on the ground field, of two indecomposable $\uq$-modules.
The extendable modules, together with the $R$-matrix of
$\mathrm{u}_q$ and
the action of the Auslander-Reiten transpose (see \cite{1}) on
the category of modules, complete the proof which was not achieved in \cite{2}.  
As a consequence the tensor
product commutes, despite the non quasi cocommutativity of $\uq.$
Moreover we obtain explicit isomorphisms between $M\otimes N$ and
$N\otimes M$ for any two $\uq$-modules and we can observe that these
canonical isomorphisms  have the properties of morphisms in a braided category 
(see \cite{5}), except of course that they are not natural.\\
We also consider tensor products of simple modules over the entire $\mathrm{u}_q$. 
The crucial observation is that extendable non-projective
$\uq$-modules are the simple modules on $\mathrm{u}_q$. 
A connection between the decomposition formulas over $\uq$ and $\mathrm{u}_q$ 
is established. 
We thus derive  formulas previously obtained by Reshetikhin and Turaev
in \cite{8} for the tensor product of simple 
$\mathrm{u}_q$-modules in a new way. The proof we obtain is new and  entirely based
on basic properties of extendable modules.

Furthermore we establish a totally different proof of the
decomposition formula for $\uq$-modules which actually includes the
three situations $\uq$, the universal enveloping algebra $U(\frak{sl}_2)$ of
$\frak{sl}_2$ and the quantum universal enveloping algebra $U_q(\frak{sl}_2)$ of
$\frak{sl}_2$  when  $q$ is not a root of unity. The proof consists in a fairly simple
axiomatisation on the Grothendieck ring of these Hopf algebras.
\section{The Hopf algebras $\mathrm{u}_q(\frak{sl}_2(\Bbb C))$ and u$_q^+$.}
We recall definitions and known facts about the above algebras, choosing Kassel's (see \cite{5})
presentation of $\mathrm{u}_q$, originally from Lusztig (see \cite{7}). 
Let $q$ be a primitive $n$-th root of unity in $\C$, $n$ different from $4$, and set\\ 
$$d=\left\{\begin{tabular}{ll}
n&if n is odd\\
n/2& if n is even\\
\end{tabular}\right.$$ 
\begin{defi}The Hopf algebra $\mathrm{u}_q(\sld)$ 
is defined over $\C$ by the generators {$E, F, K$} and the relations :
$$E^d=F^d=0 ,\,\ \  K^d=1 ,\,\ \   KE=q^2EK ,\,\ \ KF=q^{-2}FK$$
 $$\mathrm{and}\ \ EF-FE=\frac{K-K^{-1}}{q-q^{-1}}.$$
\end{defi}
It admits a Poincaré-Birkhoff-Witt type basis in the set $\{
E^iK^jF^l\}$ for $0\leq i,j,l\ \leq d-1$ (see \cite{5}).\\
The coalgebra structure is given on the generators as follows~:\\
 the comultiplication $\Delta \ :\ \mathrm{u}_q\longrightarrow \mathrm{u}_q\otimes
\mathrm{u}_q$\,  is defined by  
$$ \begin{array}{lcl}
\Delta ( E) & = & 1\otimes E+E\otimes K\\
\Delta ( F) & = & K^{-1}\otimes F+F\otimes 1 \\ 
\Delta (K) & = & K\otimes K,
 \end{array}$$
the counit \,  
$\epsilon\ :\ \mathrm{u}_q \longrightarrow\ k$ \,  by
$$\epsilon (E)=\epsilon (F)=0\ \ \epsilon (K)=1$$
and the antipode, 
$S:\mathrm{u}_q\longrightarrow \mathrm{u}_q$, is given by 
$$S(E)=-EK^{-1} \, ,\, 
S(F)=-KF \, ,\, 
S(K)=K^{-1}$$
We have the following formulas for the comultiplication :\\
$$\Delta(E)^r=\sum_{k=0}^{j}q^{-k(r-k)}
   \begin{bmatrix}
r-k \\ r
   \end{bmatrix}_qE^k \otimes K^kE^{r-k}\ \ \mathrm{and}$$
$$\Delta(F)^r=\sum_{k=0}^{r}q^{k(r-k)}
   \begin{bmatrix}
r-k \\ r
   \end{bmatrix}_qF^kK^{-(r-k)} \otimes F^{r-k}$$
Where 
   $\begin{bmatrix}
x \\ y
   \end{bmatrix}
= \frac{[y]!}{[x]![y-x]!}$ with $[x]!=[1][2]\ldots[x] $ and 
$[x]=\frac{q^x-q^{-x}}{q-q^{-1}}.$ A formula which calculates the commutators $[E^m,F^m]$
with $m\in\{0,\ldots, d-1\}$ will be needed (see \cite{5}):
$$E^mF^m=\sum_{h=0}^{m}c_hF^{m-h}E^{m-h}\prod_{j=0}^{h-1}\frac{Kq^{-j}-K^{-1}q^j}{q-q^{-1}}$$
where $c_h$ is a nonzero coefficient.
It is well known that this Hopf algebra is quasi-triangular (see \cite{5}, \cite{6}). Its  $R$-matrix
has the following expression     :
$$R=\frac{1}{d}\sum_{0\leq i,j,k\leq d-1}\frac{(q-q^{-1})^k}{[k]!}q^{k(k-1)/2+2k(i-j)-2ij}
E^kK^i\otimes F^kK^j$$
   \begin{rk} \emph{1) \ Hopf algebras have the property that the tensor
product over the ground field of two left modules is still a left
module. Indeed, for a Hopf algebra $H$,  restricting the natural action
of $H\otimes H$ to $H$ through the comultiplication $\Delta$ yields a left
$H$-module structure.\\
2)  Recall that the $R$-matrix satisfies in particular the relation 
$\Delta^{op}=R\Delta R^{-1}$,
 where $\Delta^{op}=\tau\Delta$ and 
$\tau$ is the flip, $\tau(a\otimes b)=b\otimes a$ for $a,b\in \mathrm{u}_q.$
This relation is equivalent to the existence of a family of natural isomorphisms
between $M\otimes N$ and $N\otimes M$ for any $\mathrm{u}_q$-modules $M$ and $N$.
The isomorphisms are given by the action of $\tau R.$}
\end{rk}
The upper triangular sub-algebra of $\mathrm{u}_q$ generated by $E$ and $K$ is a
sub-Hopf algebra, denoted by u$_q^+$ ; indeed \\
$\Delta (\hbox{u}_q^+)\subset \hbox{u}_q^+\otimes\hbox{u}_q^+$ \\
$S(\hbox{u}_q^+)\subset \hbox{u}_q^+.$\\ 
The dimension over $k$ of u$_q^+$ is  $d^2.$ The set  $\{E^iK^j\}_{0\leq
i,j\leq d-1}$ is a basis of $\uq$ (see \cite{7}).
   \begin{rk}
\emph{In \cite{2} it has been shown that $\uq$ is isomorphic to a quotient
of a path algebra endowed with a Hopf algebra structure. It is our
reference for the following remarks as well as for the representations of
$\uq.$\\
1) As an associative algebra $\uq$ is uniserial, meaning
that each indecomposable module on $\uq$ has a unique decomposition
series. As a consequence $\uq$ is of finite representation type.\\
2) The Jacobson radical of $\uq$ is generated by $E.$} 
\end{rk}
We have the following proposition :
\begin{pro}
If $q$ is an $n$-th root of unity with $n\neq 4$, the Hopf algebra $\uq$ is not quasi-cocommutative.
\end{pro}
{\bf Proof :}\,  Suppose there exists an invertible
element $R\in
\uq\otimes\uq $ such that $\Delta^{op}=R\Delta R^{-1}$. Then $R$ is of
the form
$R=\sum_{0\leq i,j,k,l\leq d-1}a_{i,j,k,l}E^iK^j\otimes E^kK^l$ \, where the 
$a_{i,j,k,l}$ belong to $k$.
We have in  particular $\Delta^{op}(K)=R\Delta(K)R^{-1}$, i.e. $K\otimes KR=RK\otimes K$, 
implying that  $a_{i,j,k,l}q^{2(i+k)}=a_{i,j,k,l}.$ Hence the expression of $R$ 
must reduce to\\ 
$R=\sum_{0\leq i,j,l\leq d-1}a_{i,j,d-i,l}
E^iK^j\otimes E^{n-i}K^l\ +\ \sum_{0\leq j,l\leq
n-1}a_{0,j,0,l}K^j\otimes K^l$.\\
 In order to show that the coefficients
$a_{o,j,o,l}=a_{j,l}$ are $0$ we use the identity $\Delta^{op}(E)R=R\Delta (E)$ and
obtain the relations\\
 $a_{j,l}=q^{2j}a_{j,l-1}$ and 
$a_{j,l}=q^{-2l}a_{j-1,l}$  whenever they make sense. As a consequence, $a_{j,o}=a_{0,l}=a_{0,0}$, 
implying that $a_{1,1}=q^2a_{1,0}=q^{-2}a_{0,1}$,  hence $a_{1,1}=0$ and
$a_{j,o}=a_{0,l}=0.$ We infer $a_{i,j}=0$ \, for all $ 0\leq
j,l\leq d-1$, and $R$ is therefore reduced to 
$R=\sum_{0\leq i,j,l\leq d-1}a_{i,j,d-i,l}
E^iK^j\otimes E^{n-i}K^l.$\\
Finally, we note that $\Delta^{op}(E^{d-1})$ must be different from
zero, and then develop the expression $R\Delta(E^{d-1}).$\\ 
Writing
$\Delta(E^{d-1})=\sum_{0\leq x,y,z\leq d-1}b_{x,y,z}E^xK^y\otimes E^{d-1-x}K^z\ $\, 
with\,  $\ \ b_{x,y,z}\in \C$, we obtain $R\Delta (E^{d-1}) =
 \sum_{i,j,l,x,y,z}c_{i,j,l,x,y,z}E^{i+x}K^{j+y}\otimes E^{2d-i-x-1}K^{l+z}.$
Since either $i+x\geq d$ or $2d-i-x-1\geq d$, we necessarily have
$R\Delta (E^{d-1})=0.$ We thus arrive to the contradiction
 $R\Delta(E^{d-1})R^{-1}=0$ and
$R\Delta(E^{d-1})R^{-1}=\Delta^{op}(E^{d-1})\not=0.$  $\ \square$ 
\begin{rk}\emph{The case $n=4$ yields a quasi-cocommutative Hopf
algebra (see \cite{2}). An alternative proof of {Proposition 2.1}
is provided in \cite{2} using the presentation of  $\uq$
by quiver and relations .}
\end{rk}
\subsection{Modules.}
The isomorphism classes of the modules described below
constitute the complete list of isomorphism classes of indecomposable
$\uq$-modules;  they are all non-isomorphic. 
To each couple $(i,u)$, where $i\in\Z/d\Z$ and $0\leq u\leq d-1$,
corresponds a $\uq$-module, denoted by $M_i^u$, of dimension $u+1.$ It
admits  a basis 
$\{e_i^0,e_i^1,\ldots , e_i^u\}$  over $\C$ such that the action of $\uq$ on the
basis vectors is given by
$$\left\{
\begin{array}{lcl}
Ke_i^j & = & q^{2(i+j)}e_i^j\\
Ee_i^j &  = & e_i^{j+1}\ \  \mathrm{for} \ 0\leq j\leq u-1\\
Ee_i^u & = & 0
\end{array}
\right.$$
Note that $e_i^0$ is a generator of $M_i^u$ over $\uq$.
The indecomposable projective modules are those
of dimension $d-1$, and we
denote them by $P_i=M_i^{d-1}$. The simple modules are the one-dimensional
modules, and we denote them by $S_i=M_i^0.$\\
{\bf Notations} : The length of a vector $v$ belonging to a $\uq$-module
is an integer $0\leq m\leq d-1$, minimal for the property $E^{m+1}v=0.$
In particular, the length of a basis vector of the type $e_i^j$ 
is $u-j.$ For $r,s\in \Z$ let E$(r/s)$ be the entire part of $r/s.$ 
\section{Axiomatisation of the tensor product of modules.}
The tensor product of modules on $\uq$ has decomposition
formulas which are similar to those for the universal enveloping algebra of
$\sld$, and for the quantum universal enveloping algebra of $\sld$ when
$q$ is not a root of unity. The following axiomatisation unifies the
proofs of these formulas leaving behind the concrete decomposition.
\begin{rk} 
\emph{Recall that the Grothendieck group of a ring $\Lambda$,
denoted by $K(\Lambda)$, is the quotient of the free abelian group
with basis the
isomorphism classes $[X]$ of modules $X$ on $\Lambda$ by the subgroup
generated by elements \\
$[X_2]-[X_1]-[X_3]$ provided by each split exact
sequence $X_1\rightarrow X_2\rightarrow X_3$ of \\
$\Lambda$-modules.   
Moreover if $\Lambda$ is a Hopf algebra, the free abelian group is endowed
with a ring structure through the tensor product of modules. 
The functor induced by tensoring over the ground field is flat,
implying that the subgroup above is an ideal, and hence 
 the quotient $K(\Lambda)$ is still a ring.\\
If $\Lambda$ is a finite dimensional algebra, its Grothendieck group
is a free abelian group with basis given by the isomorphism classes of
indecomposable modules.}
\end{rk}
Let $I$ be the set $\{0\}$ or $\Z/d\Z.$ To $m$ belonging to 
$\overline{\N}-\{0\},$ where $\overline{\N}=\N\cup\infty,$ put $J_m$ to be the set
$\{0,\ldots , m-1\}$ if $m\in\N$ and $J_m=\N$ if $m=\infty.$
 Consider  the free
commutative group generated by the elements [i,u], where $(i,u)$
 belong to $ I\times J_m.$ Suppose now that this group
is equipped with an extra
multiplicative structure, making it into a ring.  Denote by $\oplus$ the
addition law and $\otimes$ the multiplication law. We need to put $[i,u]=0$ if $u<0.$
We have the  proposition : 
\begin{pro}
Assume the  relations below hold
and are symmetric with respect to $\otimes$ :\\
$[i,0]\otimes [j,0]=[i+j,0]$, 
$[0,1]\otimes [j,v]=[j,v+1]\oplus[j+1,v-1]$ \, for
$0\leq v\leq m-2$ \, , \\
and $[0,1]\otimes [j,m-1]=[j,m-1]\oplus[j+1,m-1]$ where $(i,j)\in I\times J_m$ and 
$u,v\in J_m.$\\
Then the following decomposition formulas are true:
   \begin{enumerate}
      \item $[i,u]\otimes [j,v]=\oplus_{l=0}^{min(u,v)}[i+j+l,u+v-2l]$ \, for
$u+v\leq m-1$ 
      \item $[i,u]\otimes [j,v]=\oplus_{l=0}^e[i+j+l,m-1]\oplus 
\oplus_{l=e+1}^{min(u,v)}[i+j+l,u+v-2l]$ \, for $u+v\geq m-1$ \, where 
$e=u+v-(m-1)$
\end{enumerate}
   \end{pro}
{\bf Proof : }\, We proceed by double induction. First we prove that\\ 
$[i,u]\otimes [j,0]=[i+j,u]$ for all $i,j\in I$ and $u\leq m-1$ by
induction on $u.$ By
assumption it is true for $u=0$. Suppose it is valid up to a rank $0<u<m-1$ and
let's show it for $u+1.$ For this purpose we look at 
$[0,1]\otimes [i,u]\otimes [j,0].$ Developing the left and right side
respectively we obtain the equality :\\
$([i,u+1]\otimes [j,0])\oplus
([i+1,u-1]\otimes[j,0])=[0,1]\otimes[i+j,u]$\\
that is $([i,u+1]\otimes
[j,0])\oplus[i+j+1,u-1]=[i+j,u+1]\oplus[i+j+1,u-1]$,\\
and as a consequence $[i,u+1]\otimes [j,0]=[i+j,u+1].$\\
Next, we take an arbitrary $u$, and  show the formulas by
induction on $v.$ Suppose they hold up to a rank $v\geq 1$, then we have two
situations to consider, either $u+v+1\leq m-1$ or $u+v+1\geq m-1.$
Developing $[0,1]\otimes [i,u]\otimes [j,v]$ on the left and right
hand side respectively easily solves the first case. For the second
more care is needed. Set $e=u+v-(m-1)$ and let $\oplus_{l=0}^a [x_l,y_l]=0$ if  $a\leq 0 $
with $(x_l,y_l\in I\times J_m).$ 
We proceed as before by developing the left and right sides of
$[i,u]\otimes [j,v]\otimes [0,1]$ and thus obtaining the
equality\\
 $(\oplus_{l=0}^{e}[i+j+l,m-1]\oplus\
\oplus_{l=e+1}^{min(u,v)}[i+j+l,u+v-2l])\otimes
[0,1]=$\\
$[i,u]\otimes([j,v+1]\oplus[j+1,v-1]).$ Developing this gives us the identity\\
$(\oplus_{l=0}^{e}([i+j+l,m-1]\oplus [i+j+l+1,m-1])
\oplus $\\
$ \oplus_{l=e+1}^{min(u,v)}([i+j+l,u+v-2l+1]\oplus
[i+j+l+1,u+v-2l-1])$ = $[i,u]\otimes [j,v+1]\oplus \oplus_{l=0}^{e-1}[i+j+1+l,m-1]\oplus\
\oplus_{l=e}^{min(u,v-1)}[i+j+1+l,u+v-1-2l]).$\\
Therefore
$[i,u]\otimes[j,v+1]=\oplus_{l=0}^{e+1}[i+j+l,m-1]\oplus
\oplus_{l=e+2}^{min(u,v+1)}[i+j+l,u+v+1-2l].\ \ \square$
   \begin{rk} 
\emph{The Grothendieck ring of the Hopf algebra $\uq$
corresponds to $I=\Z/d\Z$ and $m=d$ where we replace the formal writing $[i,u]$
by the isomorphism class of the indecomposable module $[M_i^u].$ This observation
leads us to the next result.}
   \end{rk}
   \begin{theo}
Let $M_i^u$ and $M_j^v$ be indecomposable $\uq$-modules for\\
$i,j\in\Z/d\Z$ \, and\,  $0\ \leq u,v\ \leq d-1$.  There are isomorphisms :
   \begin{enumerate} 
     \item 
If $u+v\ \leq\ d-1$
$$M_i^u\otimes M_j^v\cong \bigoplus_{l=0}^{min(u,v)}M_{i+j+l}^{u+v-2l}$$
     \item
If $u+v\ \geq\ d-1$, set $e=u+v-(d-1)$, then
$$M_i^u\otimes M_j^v\cong \bigoplus_{l=0}^eP_{i+j+l}\oplus\bigoplus_{l=e+1}^{min(u,v)}M_{i+j+l}^{u+v-2l}$$ 
   \end{enumerate}
    \end{theo} 
{\bf Proof :}\, In view of the previous remark we can apply
the   proposition. We need to  check
that $S_i\otimes S_j\cong S_{i+j}$ and  $M_i^1\otimes S_j\cong M_{i+j}^1$ as well as $M_i^1\otimes M_j^{v}\cong
M_j^{v}\otimes M_i^1\cong M_{i+j}^{v+1}\oplus M_{i+j+1}^{v-1}$ and finally
that $M_i^1\otimes P_j\cong P_j\otimes M_i^1\cong P_{i+j}\oplus P_{i+j+1}.$\\
The first two isomorphisms are simply given by letting $e_i^0\otimes e_j^0$
go to a non zero multiple of $e_{i+j}^0.$\\
To prove the third assertion (we treat the case $M_i^1\otimes
   M_j^{v}\cong  M_{i+j}^{v+1}\oplus M_{i+j+1}^{v-1}$) we need to ensure 
that in $M_i^1\otimes M_j^v$ we have two vectors $w_1$ and $w_2$, 
 of lengths $v+1$ and  $v-1$ respectively, and whose $K-$eigenvalues are
respectively $q^{2(i+j)}$ and $q^{2(i+j+1)}.$ Indeed this implies the
existence of $M_{i+j}^{v+1}$ and $M_{i+j+1}^{v-1}$ as submodules of
$M_{i}^1\otimes M_j^v$, as well as their sum which is necessarily
direct. For dimension reasons we therefore obtain the required
isomorphism.\\
Let us make explicit the vectors $w_1$ and $w_2$. For $w_1$ we simply take
   $e_i^0\otimes e_j^0.$ What needs to be checked is that $E^{v+1}e_i^0\otimes
   e_j^0\neq 0$ (note that $E^{v+2}e_i^0\otimes e_j^0$ is necessarily
   equal to 0). Using the comultiplication formulas given in section 2 we
   find that $E^{v+1}e_i^0\otimes e_j^0=q^{-v}
\begin{bmatrix}
v+1\\
1
\end{bmatrix}_qq^{2(v+j)}e_i^0\otimes
   e_j^v$ ; this is equal to
   $\frac{q^{v+1}-q^{-v-1}}{q-q^{-1}}e_i^1\otimes e_j^v$ which is not
   equal to $0$
   since we are in the case $v\leq d-1.$ To determine $w_2$ we need to
   make two computations~: first, let $a,b$ belong to $k$, then we have
   $E^{v-1}(ae_i^1\otimes e_j^0 + e_i^0\otimes e_j^1)= be_i^0\otimes
   e_j^v +(a+b(q^{2j+v-2})\frac{q^{v-1}-q^{-(v-1)}}{q-q^{-1}})e_i^1\otimes
   e_j^{v-1}$, which is non-zero whenever $a$ and $b$ are both different from
   zero. Next, we compute $E^v(ae_i^1\otimes e_j^0+be_i^0\otimes
   e_j^1)$ and find it to be equal to
$ae_i^1\otimes e_j^v
   + q^{2j+v+1}\frac{q^v-q^{-v}}{q-q^{-1}}be_i^1\otimes e_j^v.$ In
   view of these computations, we set $w_2=
   ae_i^1\otimes e_j^0 +be_i^0\otimes e_j^1$, with 
$a=-q^{2(j+v)}+1+q^{2j+1}$ and $b=q-q^{-1}$, and hence obtain a vector
   satisfying the desired conditions.
$\square$
  \begin{rk}
\emph{ We will see that the theorem can be obtained in a
totally different way, by means of extendable $\uq$-modules.}
\end{rk}
Next we consider two different cases where our axiomatisation applies.
  \begin{pro}
   Taking $I={0}$ and $m=\infty$ leads to Clebsch-Gordan
formulas for $U(\sld)$ and $U_q(\sld)$ when $q$ is not a root of
unity. 
   \end{pro}
{\bf Proof :}\,1)  Recall the irreducible representations of $U(\sld)$. To each
integer $n$ corresponds a simple $U(\sld)$-module $V(n)$ of dimension
$n+1.$ It admits a basis $  \{v_0,\ldots,v_n\}$ over $k$ such that the action of
$U(\sld)$ is given by
$$\left\{
  \begin{array}{lcl}
Xv_i & = & (n-i+1)v_{i-1}\\
Yv_i & = & (i+1)v_{i+1}\\
Hv_i & = & (n-2i)v_i
  \end{array} \mathrm{where}\  v_i=0\ \mathrm{for}\ i\not\in\{0,\ldots,n\}  
\right.$$
and we have the Clebsch-Gordan formula for the decomposition of the
  tensor product of two such modules :
$V(n)\otimes V(m)\cong \oplus_{l=0}^{min(n,m)}V(n+m-2l).$ In view of
  the preceding results, this formula can be obtained by checking the
  following isomorphisms of $U(\sld)$ :
$V(0)\otimes V(0)\cong V(0)$ and\\
$V(1)\otimes V(n)\cong V(n+1)\oplus V(n-1)$ for $n\geq 1.$\\
The first is trivial, the second is obtained by giving an
explicit decomposition as it was done for $\uq$. Indeed, let $\{v_0,v_1\}$
and $\{ v'_0,\ldots, v'_n\}$ be the basis of $V(1)$ and $V(n)$
respectively. Then the vectors  $v_0\otimes  v'_0$ and $v_0\otimes
v'_1- mv_1\otimes v'_0$ are generators of the modules $V(n+1)$ and
$V(n-1)$ respectively. Their sum is a direct sum and comparing the
dimensions  leads to the desired isomorphism.\\
2) The case of $U_q(\sld)$ when $q$ is not a root of unity is
similar. Let $\epsilon =\pm 1$. To each integer $n$ correspond two
modules $V_{1,n}$ and $V_{-1,n}$ who admit bases
  $\{v_{\epsilon,0},\ v_{\epsilon,1},\ldots,v_{\epsilon,n-1}\}$ such that
  the action of $U_q(\sld)$ is given by
$$\left\{
  \begin{array}{lcl}
Ev_{\epsilon,i} & = & \epsilon[n-i+1]v_{\epsilon,i-1}\\
Fv_{\epsilon,i} & = & \epsilon[i+1]v_{\epsilon,i+1}\\
Kv_{\epsilon,i} & = & \epsilon q^{n-2i}v_{\epsilon,i}.
  \end{array}\mathrm{where}\  v_{\epsilon,i}=0\ \mathrm{for}\ i\not\in\{0,\ldots,n\}
\right.$$  
The Clebsch-Gordan formula is : $V_{\epsilon,n}\otimes
  V_{\epsilon',m}\cong
  \oplus_{l=0}^{min(n,m)}V_{\epsilon\epsilon',n+m-2l}.$
One easily reduces to the case of modules of type $V_{1,n}$ and as in
  the former situations the isomorphism between $V_{1,1}\otimes V_{1,n}$
  and $V_{1,n+1}\oplus V_{1,n-1}$ for $n\geq 1$ is guaranteed by the two
  vectors $v_0\otimes v'_0$ and $v_0\otimes v'_1-[m]q^{-m}v_1\otimes
  v'_0$ (we assume that the vectors $v_i$ and $v'_j$ form bases for $V_{1,1}$ and
  $V_{1,n}$ respectively).
$\square$
 \begin{rk} 
  \emph{Considering the simple $\uq$-modules, we can observe that they
  form a multiplicative group for the tensor product,
  isomorphic to the cyclic group of order  $n.$ Actually, the
  isomorphism classes of  
  simple modules over a basic and split Hopf algebra always provide a group 
  (see for instance \cite{4}). Now
  this group acts
  on the category of $\uq$-modules via the tensor product and it is
  interesting to note that 
  the action of the generator $S_1$ on an indecomposable module yields
  the dual transpose (see \cite{1}).}
\end{rk} 
%
\section{Extendable modules.}
It is obvious that a $\uq$-module is not in general  issued from a
 $\mathrm{u}_q$-module, in the sense that it is not obtained by
  restricting the action of $\mathrm{u}_q$ to $\uq$. Nevertheless we
  can consider the subfamily of $\uq$-modules on which indeed there
  exists an action of $\mathrm{u}_q$ such that the original action of
  $\uq$ is respected. We call those modules extendable. They have the
  property that the $R$-matrix of $\mathrm{u}_q$ provides
  isomorphisms making the tensor product of two such
  modules commutative. Restricting our study to this family gives some
  information on the decomposition of $\uq$-modules, as well as on
  simple $\mathrm{u}_q$-modules. We need the following notation :\\
{\bf Notation} : Let $u\in \N$, then $\overline u$ is the representative
element of the class of $u$ modulo $d$ contained in the set $\{0,\ldots ,d-1\}.$
  \begin{theo}
The extendable indecomposable modules are :
  \begin{enumerate}
\item
The indecomposable modules of type $M_i^{\overline{-2i}}$ for
   $0\leq i\leq d-1$.   These modules extend in a
 unique way and provide all the simple $\mathrm{u}_q$-modules.
\item 
The projective indecomposable modules  $P_i$ \, for $0\leq i\leq d-1$. 
 These
  modules extend in two non-isomorphic ways, except $P_{\frac{d+1}{2}}$ 
when $d$ is odd.
   \end{enumerate}
  \end{theo}
{\bf Proof :}\, 
We proceed in the following way : First we consider an arbitrary
indecomposable $\uq$-module, and we try to define an action of $F\in
\mathrm{u}_q$ on its basis elements, such that the original action of
$\uq$ is preserved, and the algebra structure of $\mathrm{u}_q$ is
respected. We thus infer the necessary conditions for an indecomposable 
module to be extendable.\\
Consider a module  $M_i^u$ with $i\in\Z/d\Z\ $ and
$\ 0\leq u\leq d-1$.
It is generated over $\uq$ by the element $e_i^0$, and the set $\{e_i^j\}_{0\leq j\leq u}$ 
is a basis over $k$. 
The action of $\uq$ is given by $Ee_i^j=e_i^{j+1}$
for $0\leq j\leq u-1\ $,  $Ee_i^u=0$ and $Ke_i^j=q^{2(i+j)}e_i^j$.\\
Suppose we have an action of $F$ given by $Fe_i^j=\sum_{o\leq h\leq u}\lambda_{i,j}^he_i^h$ 
where\\
$\lambda_{i,j}^h\in \C$.  The relation $KF=q^{-2}FK$ implies
$KFe_i^j=\sum_{0\leq h\leq
u}\lambda_{i,j}^hq^{2(i+h)}e_i^h\\=q^{-2}FKe_i^j
=q^{-2}q^{2(i+j)}\sum_{0\leq h\leq u}\lambda_{i,j}^he_i^h.$\\
It follows that $\lambda_{i,j}^hq^{2(i+h)}=q^{2(i+j-1)}\lambda_{i,j}^h$
and therefore :\\
$Fe_i^j=\lambda_{i,j}^{j-1}e_i^{j-1}=\lambda_i^{j-1}e_i^{j-1}$
for all $1\leq j\leq u$ and $Fe_i^0=\lambda_i^{d-1}e_i^{d-1}$,\\
where $\lambda_i^{d-1} =0$ if $u\leq d-2.$
Since $EF-FE=\frac{K-K^{-1}}{q-q^{-1}},$ 
 we must have the following~:
$$\lambda_i^0=-[2i]+\lambda_i^{d-1}$$ 
 We next proceed by induction and obtain
$$\lambda_j=\sum_{0\leq h\leq j}-[2(i+h)]+\lambda_i^{d-1}.$$
The remaining relations are now $Ee_i^u=0$ and $F^d=0$. From the first one
we deduce :
$$(EF-FE)e_i^u=E\lambda_i^{u-1}e_i^{u-1}=\lambda_i^{u-1}e_i^u
=\frac{K-K^{-1}}{q-q^{-1}}e_i^{u}=\frac{q^{2(i+u)}-q^{-2(i+u)}}{q-q^{-1}}e_i^{u}.$$
On the other hand, $\lambda_{u-1}=\sum_{0\leq h\leq
u-1}-[2(i+h)]+\lambda_i^{d-1}.$
The equality is automatically realized when dealing with a projective module.
Otherwise, that is when $u\leq d-2$, we need\\
 
\begin{tabular}{ll}
&$\sum_{0\leq h\leq
u}\frac{q^{-2(i+h)}-q^{2(i+h)}}{q-q^{-1}}$ \\ \\
=&$\frac{-q^{2i}}{q-q^{-1}}(\frac{1-q^{2(u+1)}}{1-q^2})+
\frac{q^{-2i}}{q-q^{-1}}(\frac{1-q^{-2(u+1)}}{1-q^{-2}})$\\
\\
=&$\frac{q^{2i}(1-q^{2(u+1)}+q^{-4i+2}-q^{-2(u+1)+2-4i})}{(q-q^{-1})(q^{2}-1)}$\\
\\
=&$\frac{q^{2i}(1-q^{2(u+1)})(1-q^{-2u-4i})}{(q-q^{-1})(q^{2}-1)}=0$\\ \\
\end{tabular}\\
The equality is true when $2(u+1)=0\ \hbox{mod}\ n$
and $2u=-4i\ \hbox{mod}\ n$.
For  $n$ odd the first case is never realized, and for n even it
corresponds to the projective modules. Otherwise we need the condition 
$u=-2i\ \hbox{ mod}\ d$. \\
The last condition on the $\lambda_i^j$ coming from $F^d=0$ is $\lambda
_{\frac{d+1}{2}}^{d-1}=0$ for $d$ odd.  
Hence the indecomposable modules for which the action of u$_q^+$
extends to $\mathrm{u}_q$ are the
projectives and the modules of the type $M_i^{-\overline{2i}}$. For $i\in\{0,\ldots,d-1\}$,
 it is easy to check that the modules obtained on $\mathrm{u}_q$ from
 the modules $M_i^{\overline{-2i}}$ are simple, and we thus obtain all the simple modules
 on $\mathrm{u}_q$ up to isomorphism (the list of simple $\mathrm{u}_q$-modules is given
 in \cite{5}).$\square$
  \begin{rk}\emph{ The projective modules are examples of modules extendable
to $\mathrm{u}_q$-modules in two non-isomorphic ways. We are therefore allowed to imagine the case
of an extendable module whose indecomposable components are not extendable. This 
turns out to be impossible.}
\end{rk}
 \begin{pro} A $\uq$-module is extendable if and only if it is a direct
sum of indecomposable extendable modules.
  \end{pro}
{\bf Proof :}\, 
Let $X$ be an arbitrary $\uq$-module, decomposable into
$M_i^u\oplus(\bigoplus_{l\in L\ v\in V} M_l^v)$, where $L$ is a finite set.
 We examine  the possible
actions of $F$ on the basis $\{e_i^j\}$ of $M_i^u.$ Using
a simple induction and the relation
$EF-FE=\frac{K-K^{-1}}{q-q^{-1}}$, we find that an action must be of
the form :\\ $F.e_i^j=\lambda_i^{j-1}e_i^{j-1}\ +$ linear
combination of $\{e_l^k\}_{j\leq k}.$\\ The action of $E$ on $e_i^u$
given by  $Ee_i^u=0$ requires that\\ 
   \begin{tabular}{ll}
$(EF-FE)e_i^u$&=$\frac{q^{2(i+u)}-q^{-2(i+u)}}{q-q^{-1}}e_i^u$\\
&=$EF(e_i^u)$\\
&=$(\lambda_i^{u}e_i^u\ +$ lin.comb.$ \{e_l^k\}_{u+1\leq k})$\\
   \end{tabular}\\
This implies that
$\lambda_i^{u-1}=\frac{q^{2(i+u)}-q^{-2(i+u)}}{q-q^{-1}}$, i.e. that
$M_i^u$ is an extendable module.\ $\square$  
   \begin{rk}
  \emph{ There may be more than one way to extend a direct sum
of non projective, indecomposable, extendable modules. As an example we can
give the $\uq$-module $M_1^1\oplus S_0$ in the case $n=3$. Indeed the
possible actions of $F$ are easily found to be : $Fe_1^0=ce_0^0$,
$Fe_1^1=\lambda_1^0e_1^0$ and $Fe_0^0=0$ where $c$ belongs to
$\C$. Considering the options $c=0$ and $c\neq 0$ respectively, the
result is two non-isomorphic representations of $\mathrm{u}_q.$}
\end{rk}
   \begin{rk}
\emph{ 
- For $d$ odd there is exactly one indecomposable extendable  module per
 dimension $m$, where $1\leq m\leq d-1.$\\
-  For $d$ even there are exactly two indecomposable extendable
modules\\
 per dimension $2m+1$, where $0\leq m\leq \frac{d}{2}-1$.}
\end{rk}
The following result provides a characterisation of self-dual 
indecomposable modules in terms of extendable ones. We recall that the
dual $Hom_k(M,k)$ of a module $M$ over  a Hopf algebra $H$ over a field $k$ can be
provided with a left $H$-module structure by means of the antipode
$S$ (see \cite{10}) (we denote this left $H$-module by $^*M$)~:\\
\indent $\lambda .f(x)=f(S(\lambda )x)$ \, \, for $\lambda \in H\ ,\ f\in Hom_k(M,k)$
\, and $x\in M.$
\begin{pro}
Let $M$ be a $\uq$-module. Then the following are equivalent 
   \begin{enumerate} 
     \item 
The module $M$ is indecomposable and self-dual.
     \item
The module $M$ is indecomposable and extendable of type $M_i^u$ with\\
   $u\equiv -2i$.
\end{enumerate}
\end{pro}
{\bf Proof :\, } We consider an arbitrary indecomposable module
$M_i^u$. Let $\{ (e_i^j)^*\}$ be the dual basis of $^*M_i^u$ ; then
$(e_i^{u})^*$ is a generator of this module  and we have another basis
given by the elements $\{E^j(e_i^{u})^*\}_{0\leq j\leq u-1}.$ The action
of $K$ on $E^j(e_i^{u})^*$ is the following :\\
$KE^j(e_i^{u})^*=q^{2j}E^jK(e_i^{u})^*=q^{2(j-i-u)}E^j(e_i^{u})^*.$\\
We deduce an isomorphism between $^*M_i^u$ and $M_{n-i-u}^u.$
The explicit isomorphism is \\
$M_i^u\longrightarrow ^*M_{n-i-u}^u$\\
$e_i^j \longmapsto (-1)^jq^{j(j+2i+1)}e_{n-i-u}^*$\\
Consequently $M_i^u$ is selfdual iff $u\equiv -2i.\ \square$ 
\begin{rk}
\emph{The extendable modules provide a different proof of the 
Clebsch-Gordan-like formula for  $\uq$ stated before. 
We sketch the proof briefly.} 
\end{rk}
{\bf Proof :}\, The first step does not involve the
extendable modules (see \cite{2} where the proof is incomplete). 
It consists in showing that the tensor product of two arbitrary 
indecomposable $\uq$-modules must decompose
as follows~: $M_i^u\otimes M_j^v\cong
\oplus_{l=0}^{v}M_{i+j+l}^{x_l}$ where $u-d\leq x_l\leq u+v-d$ 
and $u+v\leq d-1$ (we request the latter condition here in order to simplify,
 and we suppose that $v\leq u$).
This is done by considering the dimension of each $K$-eigenspace and
the action of $E$ on those. Indeed the $K$-eigenvalues are
$q^{2(i+j+l)}$ with $0\leq l\leq u+v$, and the dimensions are
distributed as follows :
to $q^{2(i+j+l)}$ with $0\leq l\leq v$ corresponds a vector space of
dimension $l+1$, moreover the vector space morphism induced by $E$ between the
eigenspace of eigenvalue $q^{2(i+j+l)}$ and the one of eigenvalue
$q^{2(i+j+l+1)}$  is injective. To the same situation with $v\leq
l\leq u$ corresponds a vector space of dimension $v+1$ and the
morphism induced by $E$ is one to one. Finally, for $u\leq l\leq u+v$
the dimension is $u+v-l+1$, and $E$ induces a surjective morphism whose
kernel is one-dimensional. As a consequence the quotient by the action
 of the Jacobson radical, top($M_i^u\otimes
M_j^v$), is
$\oplus_{l=0}^{min(u,v)}S_{i+j+l}$ and we conclude by uniseriality.\\
Now in the specific case of two indecomposable and extendable
$\uq$-modules, necessarily $x_l=u+v-2l$, which is the result we want in the
general case. Indeed, to each index $i+j+l$ corresponds one and only
one extendable indecomposable module. 
Moreover the tensor product of two extendable modules
is still extendable, hence it decomposes into a direct sum of indecomposable 
extendable modules, and leaves only one choice for the value of $x_l$.
Denote by $\phi$ 
the resulting isomorphism.\\
This observation on the extendable modules immediately leads to the
solution of the general case.
Let $X$ and $Y$ be the indecomposable extendable modules
of dimension $u+1$ and $v+1$ respectively, and let $S$ be the simple module
s.t. $M_i^u\otimes M_j^v\cong S\otimes  X\otimes Y$. Then the morphism
$id\otimes \phi$ realizes the required decomposition
isomorphism.\ $\square$ 
\begin{rk}
\emph{ The $R$-matrix of $\mathrm{u}_q$ provides isomorphisms through the action of $\tau R$
between $M_i^u\otimes M_j^v$ and $M_j^v\otimes M_i^u$ when these are extendable modules.
For any simple module $S_l$, induced isomorphisms are given between 
$S_l\otimes M_i^u\otimes M_j^v$ and $S_l\otimes M_j^v\otimes M_i^u$ by 
$id_{S_l}\otimes\tau R.$ Hence explicit isomorphisms are obtained,
which make  the tensor product of any two modules lying in the orbit of 
the extendable modules under the action of the structure group commutative (see remark 3.4).
We let  Ind$\uq$ denote the set of indecomposable finite dimensional $\uq$-modules,
 and we have the following corollary.}
\end{rk}
\begin{cor}
   \begin{enumerate}
       \item When $d$ is odd, the orbit, under the action of the structure group, of the extendable
 indecomposables is \emph{Ind}$\uq$, hence isomorphisms are obtained in all cases.
       \item  When $d$ is even the orbit covers all the indecomposables whose dimension
over $k$ is odd. Hence isomorphisms are given between $M_i^u\otimes M_j^v$ 
and\ $M_j^v\otimes M_i^u$ when $u$ and $v$ are even.
   \end{enumerate}
\end{cor}
The explicit isomorphisms obtained when $d$ is odd are not natural, since $\uq$
is not quasi-cocommutative. Nevertheless they satisfy the other relations defining  
a braided module category (see \cite{5}). Denote by $c_{U,V}$ the isomorphism between 
$U\otimes V$ and $V\otimes U$, where $U,V$ are $\uq$-modules. Then we have the following :
\begin{cor}
$$c_{U,V\otimes W}=(id_V\otimes c_{U,W})(c_{U,V}\otimes id_W)$$
$$c_{U\otimes V,W}=(c_{U,W}\otimes id_V)(id_U\otimes c_{V,W})$$
$$(id_W\otimes c_{U,V})(c_{U,W}\otimes id_V)(id_U\otimes c_{V,W})=
(c_{V,W}\otimes id_U)(id_V\otimes c_{U,W})(c_{U,W}\otimes id_W)$$
\end{cor}
{\bf Proof :}\, We show the first equality, the others are obtained in a similar way.
There exist extendable modules $M_1,M_2$ and $M_3$ together with a simple module $S$
and isomorphisms :\\
$\phi_1\ :\ U\otimes V\otimes W\cong S\otimes M_1\otimes M_2\otimes M_3$\\
$\phi_2\ :\ V\otimes W\otimes U\cong S\otimes M_2\otimes M_3\otimes M_1$\\
$\phi_3\ :\ V\otimes U\otimes W\cong S\otimes M_2\otimes M_1\otimes M_3$.\\
Then\\
\begin{tabular}{ll}
$c_{U,V\otimes W}$&$=\phi_2^{-1}(id_S\otimes c_{M_1,M_2\otimes M_3})\phi_1$\\
&$=\phi_2^{-1}(id_S\otimes (id_{M_2}\otimes c_{M_1,M_3})\circ
 id_S\otimes(c_{M_1,M_2}\otimes id_{M_3}))\phi_1$\\
&$=\phi_2^{-1}(( \phi_2\circ id_V\otimes c_{U,V}\circ\phi_3^{-1})\circ
 (\phi_3\circ c_{U,V}\otimes id_{W}\circ \phi_1^{-1}))\phi_1$\\
&$=(id_V\otimes c_{U,W})((c_{U,V}\otimes id_W).\ \square$\\
\end{tabular}
   \begin{rk}
  \emph{ The underlying isomorphism of vector spaces\\ $M_i^u\otimes 
M_j^v\cong M_j^v\otimes M_j^v$
does not depend on $i$ and $j$, therefore we obtain  no new solution to the 
Yang-Baxter equation.}
   \end{rk}
\section{Tensor product of simple $\mathrm{u}_q$-modules.}
Recall that the simple $\mathrm{u}_q$-modules are obtained from indecomposable 
extendable $\uq$-modules (see proposition 4.1). We denote by $\overline M_i^u$, where
$i\in \Z/n\Z$ and $0\leq u\leq d-1$, a simple
module over $\mathrm{u}_q.$ We need to recall a family of indecomposable finite dimensional
 $\mathrm{u}_q$-modules, which are both projective and injective  (see \cite{8} and \cite{9}). 
 To begin with, take the direct sum of the
 projective indecomposable 
$\uq$-modules $P_i\oplus P_{\overline{-2i}}$, where $i\in \{0, \ldots ,\mathrm{E} ((d-1)/2)\}.$
Then we define the following action of $F$ on its basis elements, making it into a
 $\mathrm{u}_q$-module~:
$Fe_i^j=\lambda_i^{j-1}e_i^{j-1}$ and $Fe_{\overline{-2i}}^j=
e_i^{\overline{-2i}+j}+\lambda_{\overline{-2i}}^{j-1}e_{-\overline{2i}}^{j-1}$
where $j\in\{0, \ldots ,\overline{4i-1}\}$ and $Fe_{-\overline{2i}}^j=
\lambda_{\overline{-2i}}^{j-1}e_{\overline{-2i}}^{j-1}$ for $j\in\{\overline{4i-1}+1,
\ldots ,d-1\}$. We denote the resulting modules by $\tilde P_i$.
In \cite{8} Reshetikhin and Turaev give decomposition formulas for the tensor product of 
simple $\mathrm{u}_q$-modules. The proof is based on the study of indecomposable modules on $\mathrm{u}_q$ ; the Verma modules and autoinjective modules as well as exact sequences of these. These decomposition formulas are established here by a totally different approach, using the preceding results obtained on $\uq$-modules.
   \begin{theo}
 Let $\overline M_i^u$ and $\overline M_j^v$ be simple $\mathrm{u}_q$-modules 
for \, $i,j\in\Z/d\Z\ $, $0\leq u,v\leq d-1\, $ and $u+v\leq d-1.\ $ Suppose $v\leq u$.
There is an isomorphism 
$$\overline M_i^u\otimes \overline M_j^v\cong \bigoplus_{l=0}^v\overline M_{i+j+l}^{u+v-2l}. $$
   \end{theo}
{\bf Proof :\ } We simply show that there's a unique way  extending the direct sum
$X=\oplus_{l=0}^v M_{i+j+l}^{u+v-2l}$, that is by extending each module separately.\\
Recall that $M_{i+j+l}^{u+v-2l}$ is generated by $e_{i+j+l}^0$ as a $\uq$-module and 
admits the set $\{e_{i+j+l}^k\}_{0\leq k\leq u+v-2l}$ as a basis over $k$.\\
Recall also that the unique extended action of $\mathrm{u}_q$ on $M_{i+j+l}^{u+v-2l}$
is given by $Fe_{i+j+l}^k=\lambda_{i+j+l}^{k-1}e_{i+j+l}^{k-1}.$

In order to extend $X$, we study the possible actions of $F$ on the basis elements. They are 
entirely determined by the action of $F$ on the generators of each indecomposable module.
Indeed, $Fe_{i+j+l}^k=\lambda_{i+j+l}^{k-1}e_{i+j+l}^{k-1}+E^k(Fe_{i+j+l}^0).$
Let us first show that $Fe_{i+j+l}^0$ is necessarily a linear combination of elements of 
the set $\{e_{i+j+l-k}^{k-1}\}_{1\leq k\leq l-1}.$ 
Suppose $Fe_{i+j+l}^0$ is a linear combination of elements
$e_{i+j+k}^{m_k}$ with 
$0\leq k\leq v$ and $0\leq m_k\leq u+v-2k.$
Applying the identity $KF=q^{-2}FK$, we find that $m_k$ is congruent to
$l-k-1$ modulo $d.$ Therefore $m_k=l-k-1+pd$ with $p\in \Z.$ Since
$0\leq m_k\leq u+v-2k$, necessarily $p=0$ and $m_k=l-k-1.$
Consequently we can write  
$\ Fe_{i+j+l}^0=\sum_{k=1}^{l-1}a_ke_{i+j+l-k}^{k-1}\ $
with $\ a_k\in k.\ $ Using the relation
 $\displaystyle{ EF-FE=\frac{K-K^{-1}}{q-q^{-1}}} $, our previous
observation on the action of $F$ on an arbitrary basis element implies\\ 
$\displaystyle{Fe_{i+j+l}^m=\lambda_{i+j+l}^{m-1}e_{i+j+l}^{m-1}+\sum_{k=1}^{l-1}a_ke_{i+j+l-k}^{k-1+m}.}$\\
Finally, since $Ee_{i+j+l}^{u+v-2l}=0$, we must have that
$EFe_{i+j+l}^{u+v-2l}=\lambda_{i+j+l}^{u+v-2l-1}e_{i+j+l}{u+v-2l-1}.\ $ This implies\\
$\displaystyle{\sum_{k=1}^{l-1}e_{i+j+l-k}^{u+v-2l+k}=0.\ } $ But for $0\leq m\leq u+v-2(l-k)$ we have\\ 
$e_{i+j+l-k}^m\not =0,$ and since $0\leq u+v-2l+k\leq u+v-2l+2k$, we find that $a_k=0$
for $1\leq k\leq l-1.$ Hence $Fe_{i+j+l}=0$ and $Fe_{i+j+l}^m=\lambda_{i+j+l}^{m-1}
e_{i+j+l}^{l-1}\ .$ $\square$
   \begin{theo}
 Let $\overline M_i^u$ and $\overline M_j^v$ be simple $\mathrm{u}_q$-modules 
for \, $i,j\in\Z/d\Z\ $, $0\leq u,v\leq d-1\, $ and $u+v\geq d-1.\ $ 
There is an isomorphism 
$$\overline{M_i^u}\otimes \overline{M_j^v}\cong
\bigoplus_{l=0}^{\mathrm{E}(e/2)}\tilde P_{i+j+l}\oplus
\bigoplus_{l=e+1}^{min(u,v)}\overline M_{i+j+l}^{u+v-2l}$$
   \end{theo}
{\bf Proof} : \,
We can observe three cases :\\
$$\begin{array}{lll}
\left\{\begin{array}{ll}
u= & d-2i\\
v= & d-2j
\end{array}
\right. ,  &
\left\{\begin{array}{ll}
u= & 2d-2i\\
v= & d-2j
\end{array}
\right. \mathrm{and}  &
\left\{\begin{array}{ll}
u= & 2d-2i\\
v= & 2d-2j.
\end{array}
\right.   
\end{array}$$
We restrict ourselves to the first case since the only difference between these is of
elementary computational order.
We furthermore assume that min$(u,v)=v.$   
   \begin{stp}
 The tensor product decomposes in the following sum :\\
 $\overline{M_i^u}\otimes \overline{M_j^v}\cong
\overline{\oplus_{l=0}^eP_{i+j+l}}\oplus
\oplus_{l=e+1}^{v}\overline M_{i+j+l}^{u+v-2l}.$ 
   \end{stp}
{\bf Proof} :\, As in the preceding proof, $Fe_{i+j+k}^{0}$ is a linear combination of
$K$-eigenvectors with $K$-eigenvalue equal to $q^{2(i+j+k-1)}$, and
$Fe_{i+j+k}^{l}=\\
\lambda_{i+j+k}^{l-1}e_{i+j+k}^{l-1}+E^lFe_{i+j+k}^0.$
First we consider the decomposition as a $\uq$-modules decomposition
and show that the element $Fe_{i+j+l}^{0}$ is not in
$\oplus_{k=e+1}^{v}M_{i+j+k}^{u+v-2k}$ for $0\leq l\leq e.$
Indeed, for $1\leq l\leq e$  the $K$-eigenvalue of the vector $Fe_{i+j+l}^0$ is
$q^{2(i+j+l-1)}$ (note that $Fe_{i+j}^{0}=0)$, whereas for
$e+1\leq k\leq v$ and $0\leq
m_k\leq u+v-2k$ the $K$-eigenvalue for the vector $e_{i+j+k}^{m_k}$
is $q^{2(i+j+k+m_k)}$. Asking $2(i+j+k+m_k)$ to be
congruent to $2(i+j+l-1)$ modulo $n$ is equivalent to require that $k+m_k\equiv
l-1\ \mathrm{mod}\ d$. But $k+m_k\in\{ e+1,\ldots ,d-2\}$ and 
$l-1\in\{0,\ldots ,e-1\}$, therefore this congruence is impossible.  
On the other hand, a  computation similar to  that of the proof of the
preceding proposition shows that $Fe_{i+j+e+l}^0=0$ for $l=1,\ldots
, v.$ Hence the first step.
  \begin{stp}
 There exists a 
$\uq$-decomposition of $\overline{M_i^u}\otimes
\overline{M_j^v}$ such that for\\
 $0\leq k\leq \mathrm{E}(e/2)$,  the action of
$F$  on the generators $e_{i+j+k}^{0}$ of the $\uq$-modules
$P_{i+j+k}$ is zero.
   \end{stp}
{\bf Proof }\, : We show that there exists a $K$-eigenvector with
eigenvalue $q^{2(i+j+k)}$ 
(unique up to scalar multiples)
for $0\leq k\leq e$, s.t. $F$ acts on
this vector as zero. Furthermore, we show that 
for $0\leq k\leq \mathrm{E}(e/2)$, this vector is of length $d-1$
and hence generates a projective $\uq$-module.\\
The list of basis-vectors with $K$-eigenvalue equal to $q^{2(i+j+k)}$ is given by the
following set of $e+1$ vectors :\\
$\{ e_i^0\otimes e_j^k\ ,\ e_i^1\otimes e_j^{k-1},\ldots ,e_i^k\otimes
e_j^0\ , 
 e_i^{u}\otimes e_j^{d+k-u}\ ,\ e_i^{u-1}\otimes
e_j^{d+k-u+1},\ldots ,\\ e_i^{u-(e-k-1)}\otimes e_j^v\}.$\\
The action of $F$ induces a vector space morphism 
between the vector
space generated by the above vectors  and the vector space 
generated by the $e+1$ vectors of $K$-eigenvalue $q^{2(i+j+k-1)}.$
The action of $F$ is described by\\ 
$Fe_i^{m}\otimes e_j^{k-m}=q^{-2(i+m)}\lambda_j^{k-m-1}
e_i^{m}\otimes e_j^{k-m-1}\ +
\lambda_i^{m-1}e_i^{m-1}\otimes e_j^{k-m}$ and\\
$Fe_i^{u-m}\otimes e_j^{d+k-u+m}=q^{-2(u-m+i)}\lambda_j^{d+k-u+m-1}
e_i^{u-m}\otimes e_j^{d+k-u+m-1}\\ 
+ \lambda_i^{u-m-1}e_i^{u+m-1}\otimes e_j^{d+k-u+m}$,\\ 
and the corresponding matrix has the following entries:\\
$$\left\{\begin{array}{lcll}
a_{p,p} & = & q^{-2(i+p-1)}\lambda_j^{k-p}\neq 0 & \mathrm{for} \ 1\leq p\leq k-1\\
a_{p,p+1} & = & \lambda_i^{p-1}\neq 0 & \mathrm{for} \ 1\leq p\leq k-1\\
a_{p,p} &  \neq & 0 & \mathrm{for} \ k+2\leq p\leq e+1\\
a_{k+1,p} & = & 0 & \mathrm{for} \ p\neq k+2\\
0 & & \mathrm{otherwise.}
\end{array}
\right. 
$$ 
We can make the following remarks : 1) The matrix is of rank $e$ and consequently
the kernel of the morphism is one-dimensional, which gives a unique
 vector 
(up to scalar multiples),
which we denote by $v_k$, s.t. $Fv_k=0$.\, \\
2) This vector $v_k$
is a linear combination of the basis vectors $e_i^m\otimes e_j^{k-1}$, which all
appear with a non-zero coefficient. We can therefore put 
$v_k=e_i^k\otimes e_j^0 + w_k$ where $w_k$ is a linear combination of
$e_i^m\otimes e_j^{k-m}$ for $1\leq m\leq k.$\,\\
 3) The vectors $e_i^m\otimes e_j^{k-m-1}$ for 
$m\in\{0,\ldots ,k-1\}$ are all in the image of this morphism.
\\ What remains to be satisfied is that $E^{d-1}v_k\neq 0.$ For this
purpose, we write $v_k$ as above~: $v_k=\ e_i^k\otimes e_j^0+ w_k.$ 
Now there exists an integer $m$, between $0$ and $d-1$, minimal for the
property $E^{m+1}v_k=0.$ Consequently  $v_k$ generates an indecomposable
 $\uq$-module of the form
$M_{i+j+k}^m$. Since  $Fv_k=0$, this $\uq$-module  is an
extendable indecomposable $\uq$-module, and so $m=d-1$ or $m$ is
congruent to $-2(i+j+k)\ \mathrm{mod}\ d$ (thm. 4.1.). We need to exclude the second
possibility. 
Suppose that $m$ is congruent to $-2(i+j+k)$; this means that
 $m=d-2(i+j)-2k = e-2k-1$ for $\ 0\leq k\leq \mathrm{E}((e-1)/2)\ $. If $e$ is even
and $k=e/2$, then $m=d-1$, and the two situations coincide. 
Observing that $u-k\geq e-k> e-2k-1$, we compute
$E^{u-k}v_k=q^{2(u-k)}e_i^u\otimes e_j^0 \ + $( vectors linearly
independant with $e_i^u\otimes e_j^0$). Necessarily $m> u-k$, which
is a contradiction, and so $m=d-1$.\\
In $P_{i+j+k}$ with $k\in\{0,\ldots ,\mathrm{E}((e-1)/2)\}$, we put $l_k=d-2(i+j+k)+1=e-2k$,
 and we have
$Fe_{i+j+k}^{l_k}=0$ (see proof of theorem {4.1}).
   \begin{stp}
There
exists a vector $\alpha_{l_k}$ such that $F\alpha_{l_k}=e_{i+j+k}^{l_k-1}.$  
Furthermore, the $\uq$-module generated by $\alpha_{l_k}$ is isomorphic to $P_{i+j+k+l_k}.$ 
  \end{stp}
{\bf Proof} :\, We observe that $l_k\in\{\mathrm{e}((e+1)/2),\ldots ,e\}$, and since
$e_{i+j+k}^{0}$ is a linear combination of the vectors
$e_i^0\otimes e_j^{k},\ldots ,e_i^{k}\otimes e_j^{0}$, we have that
$e_{i+j+k}^{l_k-1}$ is a linear combination of the vectors
$e_i^0\otimes e_j^{l_k-1},\ldots ,e_i^{l_k-1}\otimes e_j^{0}.$ Therefore, considering the
third remark in {step 5.2.2}, there exists a vector $\alpha_{l_k}$ with $K$-eigenvalue
equal to $q^{2(i+j+k+l_k)}$ s.t. 
$F\alpha_{l_k}=e_{i+j+k}^{l_k-1}.$
We now look at the $\uq$-module generated by $\alpha_{l_k}.$ There are two things to prove :
\\
1) The module $\uq \alpha_{l_k}$ is extendable. First of all, the sum
$P_{i+j+k}+\uq \alpha_{l_k}$ of $\uq$-modules is a direct sum. In order to prove this,
we show that the vectors $e_{i+j+k}^{m}$ and $E^{m}\alpha_{l_k}$ for $m\in\{ 0,\ldots, d-1\}$
are linearly
independant. Considering their $K$-eigenvalues, this means that we must have
$E^{d-m}\alpha_{l_k}\neq a_me_{i+j+k}^{l_k-m}$ and 
$E^s\alpha_{l_k}\neq a_se_{i+j+k}^{l_k+s}$ for $m\in\{0,\ldots, l_k-1\}$
and $s\in\{l_k,\ldots, d-1\}.$ 
Indeed, if we suppose $E^{d-m}\alpha_{l_k}=a_me_k^{l_k-m}$, where $a_m$
is a nonzero coefficient, it implies $0=E^{m}E^{d-m}=a_me_k^{l_k}$, which is a contradiction.
In the same way, assume that  $E^s\alpha_{l_k}=a_se_{i+j+k}^{l_k+s}$ ; this means that 
$F^{s+1}E^s\alpha_{l_k}=0$, and therefore, in view of remark 1) in {step 5.2.2}, we have
$F^sE^s\alpha_{l_k}=b_se_{i+j+k}^{l_k}$, where $b_s\in\C.$
Applying the formula (see section 2) for the commutator
$[E^s,F^s]$, we arrive to the conclusion that $\alpha_{l_k}=c_{s}e_{i+j+k}^{l_k}$,
which is impossible. \\
2)
Now the module over $\mathrm{u}_q$ generated by $\alpha_{l_k}$ is an
extension of $P_{i+j+k}\oplus \uq \alpha_{l_k}$, hence they are both compelled to be
extendable (see {proposition 4.2}). As in the proof of {theorem 5.1},  
$\uq \alpha_{l_k}$ must be isomorphic to $ M_{i+j+k+l_k}^m$,
with $m=d-1$ or $m\equiv -2(i+j+k+l_k)\mathrm{mod}\ d$. In order to exclude the
second possibility, we suppose that $l_k=d-2(i+j+k)+1$ ; this means that 
$m\equiv -2(i+j+k+d-2(i+j+k)+1)\equiv -2(d-(i+j+k)+1)\equiv 2(i+j+k)-2\mathrm{mod}\ d.$
In this case, the vectors $E^{m}\alpha_{l_k}$ and $e_{i+j}^{d-1}$ are in
the kernel of the morphism induced by the action of $E$ on the
 vector spaces concerned. The fact that the kernel is one-dimensional gives a
contradiction and therefore $\uq\alpha_{l_k}=P_{i+j+k+l_k}$.
\begin{stp}
The $\mathrm{u}_q$-module $\overline{P_{i+j+k}\oplus P_{i+j+k+l_k}}$ is indecomposable.
\end{stp}
{\bf Proof} : Suppose it admits a non trivial decomposition
 $\overline{P_{i+j+k}\oplus P_{i+j+k+l_k}}=A\oplus B$, with $A$ and $B$ non zero.
 This implies
that as $\uq$-modules (as such we denote them by $\underline A$ and
$\underline B$)  $\underline A$ or $\underline B$ is equal to
$P_{i+j+k}$, and $\underline B$ or
$\underline A$ is equal to $P_{i+j+l_k}$ (by the Krull-Schmidt theorem).
Hence $A$ and $B$ are
extended
$\uq$-projective modules, which is excluded. $\square$

{\sc D\'epartement de Mathematiques, GTA (CNRS ESA 5030), Universit\'e
Montpellier II, case 51, 34095 MONTPELLIER, FRANCE.}\\
Email : {\tt beta@@math.univ-montp2.fr} 
\end{document}